\documentclass[amssymb,amsfonts]{amsart}
\usepackage{amssymb, latexsym}
\usepackage[T1]{fontenc}
\usepackage[intlimits,sumlimits,noDcommand,narrowiints]{kpfonts}
\usepackage[english]{babel}
\usepackage{microtype}
\usepackage{svn}
\usepackage{tensor}
\usepackage{constants}
\usepackage[pdfusetitle]{hyperref}
\usepackage[geom,sets,norm,miscops,pde]{www_math_commands}

\theoremstyle{plain}
\newtheorem{thm}{Theorem}
\newtheorem{prop}[thm]{Proposition}

\theoremstyle{definition}

\theoremstyle{remark}
\newtheorem{rmk}[thm]{Remark}

\begingroup
\catcode`\%=12
\catcode`\$=12

\gdef\GitPageFooter{$Format: v:\texttt{%h %d}; last edit: %aN on %ai.$}
\endgroup

\begin{document}
\title[$1+1$ Global Existence of TImelike Minimal Surface]{Global Existence for the Minimal Surface Equation on $\Real^{1,1}$}
\author[WWY Wong]{Willie Wai Yeung Wong}
\address{Michigan State University, East Lansing, MI, 48824}
\thanks{The author thanks the Tsinghua Sanya International Mathematics Forum, Sanya, Hainan, People's Republic of China; as well as the National Center for Theoretical Sciences, Mathematics Division, National Taiwanese University, Taipei, Taiwan, for their hospitality during the period in which this research was performed.}
\thanks{\GitPageFooter}
\email{wongwwy@member.ams.org}
\subjclass[2000]{MSC, look up on AMS}
% 35 - PDE, 35L - Hyperbolic, 35M - mixed type, 35S - psDO paraDO
% 83 - Gravity, 83C - GR, 83C20 special solutions and symmetry classes
%   83C22 EinsteinMaxwell 83C35 grav waves 83C57 black holes 83C75
%   singularities and cosmic censorship 83C60 spinor methods 83C55
%   macroscopic interaction of gravitation with matter
% 58 - Global analysis, 58J - PDE on manifolds, 58J45 hyperbolic 58J05
%   elliptic 58J40 psDO 58J47 propagation singularities

\begin{abstract}
In a 2004 paper, Lindblad demonstrated that the minimal surface equation on $\Real^{1,1}$ describing graphical time-like minimal surfaces embedded in $\Real^{1,2}$ enjoy small data global existence for compactly supported initial data, using Christodoulou's conformal method. Here we give a different, geometric proof of the same fact, which exposes more clearly the inherent null structure of the equations, and which allows us to also close the argument using relatively few derivatives and mild decay assumptions at infinity.  
\end{abstract}

\maketitle

\section{Introduction}
The equation describing graphical timelike minimal surfaces in $\Real^{1,d+1}$ can be written as the quasilinear wave equation 
\begin{equation}\label{eq:vmc}
\left( \frac{m^{\mu\nu} \phi_{,\nu}}{\sqrt{1 + m^{\sigma\tau} \phi_{,\sigma} \phi_{,\tau}}} \right)_{,\mu} = 0.
\end{equation}
In dimensions $d \geq 4$ the small data global existence follows largely from the linear decay of solutions to the wave equation which has the rate $\approx t^{-3/2}$, and is by now standard \cite{Sogge2008}. In dimensions $d = 2,3$, the linear decay of the wave equation has generic rates $\approx t^{-1/2}, t^{-1}$ respectively, which, not being integrable in time, can lead to finite-time singularity formation for even small data for generic quasilinear wave equations; see the survey article \cite{HoKlSW2014}. Equation \eqref{eq:vmc} however exhibits the null condition, which is a structural condition on the nonlinearities identified first by Klainerman \cite{Klaine1986} and Christodoulou \cite{Christ1986} in $d = 3$  and later generalized by Alinhac in $d = 2$ \cite{Alinha2001, Alinha2001a}. Using this fact Brendle \cite{Brendl2002} and Lindblad \cite{Lindbl2004} established the small data global existence for \eqref{eq:vmc} in dimensions $d = 3$ and $d = 2$ respectively. Brendle's proof followed the commuting vector field method of Klainerman \cite{Klaine1986}. Lindblad, however, gave two different proofs in his paper using respectively the commuting vector field method as well as Christodoulou's compactification method. 

In this paper we focus on the case $d = 1$. The linear wave equation exhibits no decay in $d = 1$, and hence the classical null condition cannot be used to assert that the perturbation is effectively short range, as in the case for $d \geq 2$. As a side effect this means that a direct proof of global existence for \eqref{eq:vmc} in $d = 1$ modeled after the commuting vector field method is not possible. A striking aspect of \cite{Lindbl2004} is that via the conformal compactification method, Lindblad was also able to prove the global existence of solutions to \eqref{eq:vmc} for small initial data. As Lindblad observed, the main trade-off is that for the conformal compactification method, the data must be of compact support, while in the vector field method the data is merely required to have ``sufficiently fast'' decay at infinity. The purpose of this paper is to produce an alternative proof of the $d = 1$ case allowing initial data that is not necessarily compactly supported. In the course of the discussion we will also extract some more detailed geometric information concerning the solution when the data is of compact support. 

\section{General Geometric formulation}
The result of Brendle and Lindblad are based on comparing the time-like minimal surface to a flat hyperplane by normal projection, and the scalar function $\phi$ in \eqref{eq:vmc} describes the deviation, or height, of the minimal surface as a graph over the hyperplane. The minimal surface equation can however by written intrinsically by way of the \emph{Gauss and Codazzi equations}. To quickly recall: let $M\subseteq \mathbb{R}^{1,d+1}$ be a time-like hypersurface, then the Gauss equation requires that the Riemann curvature tensor\footnote{We use the convention $\Riem[_{abc}^d]X^a Y^b Z^c = \nabla_{[X,Y]} Z^d - [\nabla_X,\nabla_Y] Z^d$ which implies $\Ricci[_{ac}] = \Riem[_{abc}^b]$.} for the induced Lorentzian metric $g$ obey\footnote{We will, throughout this paper, freely raise and lower indices using the induced metric.}
\begin{equation}\label{eq:gauss}
\Riem[_{abcd}] = k_{ac} k_{bd} - k_{ad} k_{bc} \tag{Gauss}
\end{equation}
where $k$ is the second fundamental form of the embedding of $M$. The Codazzi equation on the other hand requires
\begin{equation}\label{eq:codazzi}
\nabla_a k_{bc} - \nabla_b k_{ac} = 0 \tag{Codazzi}
\end{equation}
where $\nabla$ is the Levi-Civita connection associated to $g$. It is well known that equations \eqref{eq:gauss} and \eqref{eq:codazzi} are the only obstructions to the existence of an isometric embedding. More precisely, we have the following theorem.
\begin{thm}[Fundamental theorem of submanifolds]
Let $M$ be a simply-connected manifold, with prescribed a Lorentzian metric $g$ and a symmetric bilinear form $k$.  Suppose that connection and curvature of $g$ satisfy \eqref{eq:gauss} and \eqref{eq:codazzi}. Then there exists an isometric immersion of $M$ into Minkowski space with one higher dimension, such that $k$ is the corresponding second fundamental form. 
\end{thm}
\begin{rmk}The Riemannian version of a more general theorem is well-known. See for example Theorem 19 in Chapter 7 of \cite{Spivak1979IV}. The Lorentzian version follows \emph{mutatis mutandis} from the same proof. On the other hand, since we are considering the Cauchy problem, one can also simply ``integrate'' the second fundamental form in time once to obtain the normal vector field to $M$, and integrate it once more to get a parametrization of the hypersurface. 
\end{rmk}
The upshot of this is that, the initial value problem for \eqref{eq:vmc} can be reformulated as finding the metric $g$ and the second fundamental form $k$ from the initial data, requiring that they obey \eqref{eq:gauss} and \eqref{eq:codazzi}, 

Equations \eqref{eq:gauss} and \eqref{eq:codazzi} do not yet describe an evolution: they are underdetermined. We get a well-posed initial value problem if we combine \eqref{eq:codazzi} with the minimal surface equation $\trace_g k = 0$, which gives
\begin{equation}\label{eq:mc}
\nabla^a k_{ab} = 0. \tag{Minimality}
\end{equation}
Using that the Riemann curvature tensor is at the level of two derivatives of the metric, the initial data for the system \eqref{eq:gauss}, \eqref{eq:codazzi}, and \eqref{eq:mc} consists of the metric $g$, the connection coefficients $\Gamma$, and the second fundamental form $k$ restricted to an initial slice. The initial value problem for \eqref{eq:vmc} requires the data for $\phi$ and $\partial_t\phi$ at time $t = 0$; these determine fully the first jet of $\phi$ restricted to the initial slice. As the \eqref{eq:vmc} is second order hyperbolic, the formally provides us the $k$-jet of $\phi$ for any $k\in\Natural$ restricted to the initial slice, and a simple computation shows that the $2$-jet of $\phi$ fully determines $g$, $\Gamma$, and $k$ on the initial slice. And thus we can indeed approach the vanishing-mean-curvature evolution through the initial value problem written in terms of the Gauss-Codazzi equations. 

\section{Double null formulation when $d = 1$}
For geometric partial differential equations, it is necessary to fix a gauge (i.e\ make a coordinate choice). For hyperbolic equations it is convenient to use a double-null formulation. Since $M$ is a two-dimensional Lorentzian manifold, it is foliated by two transverse families of null curves. Equivalently, there exists two independent functions $u,v$ satisfying 
\begin{equation}\label{eq:eikonal}
g(\nabla u, \nabla u) = g(\nabla v,\nabla v) = 0
\end{equation} 
 and $g(\nabla u,\nabla v) > 0$. Note that replacing $u$ by $f(u)$ for any $f:\Real\to\Real$ with $f'\neq 0$ gives an equivalent reparametrization; we fix $u,v$ by prescribing their initial values on the initial slice, satisfying in particular that
\begin{equation}
\text{Initial slice} = \{ u - v = 0\}.
\end{equation} 
 Since $u$ and $v$ are solutions to the eikonal equation \eqref{eq:eikonal}, as is well known $L \eqdef \nabla u$ and $N \eqdef \nabla v$ are \emph{geodesic} null vector fields, at at any point $\{L,N\}$ form a null frame of the tangent space $T_pM$. 
 
The metric takes the form 
\begin{equation}\label{eq:metric}
g = (\exp \psi) (\D*{u}\otimes \D*{v} + \D*{v} \otimes \D*{u})
\end{equation}
where $\psi$ is some real valued function. This implies that 
\[ L(v) = g(\nabla u, \nabla v) = g(L,N) = N(u) = \exp (-\psi) \]
so that the coordinate derivatives are
\begin{equation}\label{eq:partials}
\partial_v = (\exp \psi) L \qquad \partial_u = (\exp \psi) N
\end{equation}
and that the inverse metric is
\begin{equation}\label{eq:imetric}
g^{-1} =  \exp(-\psi) (\partial_u \otimes \partial_v + \partial_v \otimes \partial_u) = (\exp \psi) (L \otimes N + N \otimes L). 
\end{equation}
We can also compute
\begin{equation}\label{eq:rot}
\begin{gathered}
\nabla_L N = (\exp \psi) g(L, \nabla_L N) N = - L(\psi) N,\\
\nabla_N L = (\exp \psi) g(N,\nabla_N L) L = - N(\psi) L.
\end{gathered}
\end{equation}

The assumption that $k$ is trace-free and symmetric requires that there exist scalar functions $\lambda, \nu$ such that
\begin{equation}\label{eq:sff}
k = \lambda L\otimes L + \nu N\otimes N.
\end{equation}
Then \eqref{eq:codazzi} implies
\[ (L^a N^b - N^a L^b) \nabla_{a} k_{bc} = 0 \]
and \eqref{eq:mc} becomes
\[ (L^a N^b + N^a L^b) \nabla_{a} k_{bc} = 0.\]
Taking linear combinations we get finally the \emph{null propagation equations} for the second fundamental form
\begin{equation}\label{eq:knp}
\begin{gathered}
L(\lambda) = 0, \\
N(\nu) = 0.
\end{gathered}\tag{$k$-NP}
\end{equation}
An immediate consequence of \eqref{eq:knp} is the following proposition.
\begin{prop}\label{prop:advection}
The scalar $\lambda$ is independent of $v$; the scalar $\nu$ is independent of $u$. 
\end{prop} 

The evolution of the metric is now a scalar equation for the conformal factor $\psi$. The equation \eqref{eq:gauss} implies 
\begin{equation}\label{eq:gauss2}
\Riem[_{abcd}] L^a N^b L^c N^d = \exp(-4\psi) \lambda \nu.
\end{equation}
Combined with the definition that
\[ \Riem[_{abcd}] L^a N^b L^c N^d = g( \nabla_{[L,N]} L - [\nabla_L,\nabla_N] L, N) \]
we arrive at 
\[ L(N(\psi)) + L(\psi) N(\psi) =\exp(-3\psi) \lambda \nu,\] 
which we can rewrite as
\begin{equation}\label{eq:psiw}
\partial^2_{uv} \psi = \exp(-\psi) \lambda \nu. \tag{$\psi$-wave}
\end{equation}

The equations \eqref{eq:knp} and \eqref{eq:psiw} are evolution equations on the Minkowski space $\mathbb{R}^{1,1}$ with $u,v$ the canonical null coordinates $r \pm t$, with initial data prescribed at $t = 0$. In the remainder of this paper we discuss the easy consequences from the formulation above. 

\section{Spatially compact initial data}

Suppose that the initial data for $\lambda,\nu$ is compactly supported, then for sufficiently large $u,v$ the functions $\lambda,\nu$ vanish respectively. In view of Proposition \ref{prop:advection} and \eqref{eq:gauss2}, we immediately see that
\begin{thm}
For compact initial data, the manifold $M$ is intrinsically flat outside a compact domain. 
\end{thm}
Since \eqref{eq:psiw} is now a semilinear wave equation, with the nonlinearity \emph{only supported in a fixed space-time compact region}, immediately by Cauchy stability we have
\begin{thm}
For all sufficiently small compactly supported initial data, we have global existence for the initial value problem.
\end{thm}
\begin{rmk}
The compact support assumption is only necessary for $\lambda,\nu$, and not for $\psi$!
\end{rmk}
\begin{rmk}
One can easily check that the argument closes by standard energy arguments for initial data satisfying
\begin{itemize}
\item $\lambda, \nu\in L^\infty(\{t = 0\})$ with compact support and sufficiently small norm;
\item $\psi \in H^{1/2+\epsilon}(\{t = 0\})$ and $\partial_t\psi \in L^2(\{ t = 0\})$ with sufficiently small norm. 
\end{itemize}
The $H^{1/2 + \epsilon}$ is for Sobolev embedding into $L^\infty$. The regularity outlined above is favorable compared to the classical local wellposedness level at $H^{5/2 + \epsilon}$ for the quasilinear wave equation \eqref{eq:vmc}. Recalling that the second fundamental form is roughly two derivatives of the graphical function $\phi$ and the metric is at the level of $|\partial\phi|^2$, we see that the corresponding regularity for $\lambda,\nu \in H^{1/2+\epsilon}$ and $\psi \in H^{3/2 + \epsilon}$ is higher compared to the result above. 
\end{rmk}

\section{Non-compact initial data}
We can also handle the case where $\lambda,\nu$ are not compactly supported initially, but with sufficiently strong decay. This is the main theorem of this paper. Below we let $r = u+v$ be the coordinate function on the initial slice, and $t = u-v$ the time-function on $M$; $\partial_r$ and $\partial_t$ refer to the coordinate derivatives in the $(r,t)$ coordinate system.  

\begin{thm}
Consider the initial value problem for the system \eqref{eq:knp} and \eqref{eq:psiw}, with smooth initial data $\lambda_0, \nu_0, \psi_0 = \psi|_{\{t = 0\}}$ and $\psi_1 = \partial_t \psi|_{\{t = 0\}}$. Suppose there exists $\epsilon > 0$ such that whenever the initial data satisfies 
\begin{gather*}
\norm[1]{\lambda_0} + \norm[1]{\nu_0}  \leq \epsilon,\\
\norm[\infty]{\psi_0} + \norm[1]{\partial_r\psi_0} + \norm[1]{\psi_1} < \epsilon, 
\end{gather*}
we have a unique global solution. 
\end{thm}
\begin{rmk}
As is well known, the regularity level $W^{1,1}$ is critical for wave equation in one spatial dimension. 
\end{rmk}
\begin{proof}
Using \eqref{eq:knp} and that $\lambda$ and $\nu$ propagate in different directions we immediately obtain that $\norm[L^1_{t,r}]{\lambda\nu} \approx \norm[L^1_{u,v}]{\lambda\nu} \lesssim \epsilon^2$. Using the fundamental solution for the wave equation in one dimension we get
\[ \norm[\infty]{\psi(t,\cdot)} \lesssim \epsilon + \norm[L^1({[0,t]}: L^1)]{\lambda \nu \exp(-\psi)}.\] 
This we bound by
\[ \epsilon + \exp(\norm[L^\infty({[0,t]}:L^\infty)]{\psi}) \cdot \norm[L^1({[0,t]}:L^1)]{\lambda\nu} \lesssim \epsilon + \epsilon^2 \exp(\norm[L^\infty({[0,t]}:L^\infty)]{\psi}) .\]
Then as long as $\epsilon$ is sufficiently small we have by bootstrapping that $\norm[L^\infty_{t,r}]{\psi} < 3 \epsilon$.  The theorem then follows from this a priori estimate and standard arguments.
\end{proof}

\bibliographystyle{amsalpha}
\bibliography{mixmaster.bib}

\end{document}